\documentclass[12pt]{article}
\usepackage[T2A]{fontenc}
\usepackage[english]{babel}
\usepackage{amsfonts,amssymb,amsmath,enumerate,amsthm}
\usepackage{epsfig}
\usepackage{float}

\newcommand{\RR}{\mathbb R}

\newcommand{\NN}{{\mathbb N}}

\renewcommand{\Im}{\mathop{\mathrm{Im}}}

\newtheorem{theorem}{Theorem}
\newtheorem{lemma}{Lemma}
\newtheorem{remark}{Remark}

\newtheorem{proposition}{Proposition}
\newtheorem{corollary}{Corollary}
\newtheorem{example}{Example}

\newcommand{\beq}{\begin{equation}}
\newcommand{\eeq}{\end{equation}}
\newcommand{\ba}{\begin{array}}
\newcommand{\ea}{\end{array}}
\newcommand{\bea}{\begin{eqnarray}}
\newcommand{\eea}{\end{eqnarray}}

\usepackage{graphicx}
\usepackage[usenames,dvipsnames]{color}
\usepackage{transparent}

\DeclareMathAlphabet{\mathpzc}{OT1}{pzc}{m}{it}

\usepackage{ulem}
\usepackage{cancel}

\begin{document}

\begin{center}

{\bf A holographic global uniqueness in passive imaging}
\vskip 10pt
{\it R.G. Novikov }

\vskip 10pt

{}

\end{center}

\begin{abstract}
We consider a radiation solution $\psi$ for the Helmholtz equation in an exterior region in $\RR^3$. 
We show that the restriction of $\psi$   to any ray $L$  in the exterior region is uniquely determined
by its imaginary part $\Im\psi$ on an interval of this ray.
As a corollary, the restriction of $\psi$ to any plane $X$ in the exterior region is uniquely determined   by  $\Im\psi$ on an open domain in this plane.
These results have holographic prototypes in the recent work Novikov (2024, Proc. Steklov Inst. Math. 325, 218-223).
In particular, these and known results imply a holographic type global uniqueness in passive imaging and for the 
Gelfand-Krein-Levitan inverse problem  (from boundary values of the spectral measure in the whole space) in the monochromatic case.
Some other surfaces for measurements instead of the planes $X$ are also considered.
\end{abstract}



\section{ Introduction}

We consider  the Helmholtz equation
\beq
\label{eq:1.1}
-\Delta\psi(x) = \kappa^2 \psi(x),\ \ x\in {\cal U},\ \  \kappa>0, 
\eeq
where $\Delta$ is the Laplacian in $x$, and $\cal U$ is an exterior region in $\RR^3$
that is $\cal U$ is an open  unbounded connected set  in $\RR^3$ consisting of all points  in exterior of a closed bounded regular surface $S$    (as in \cite{A}, \cite{W}).

For equation (\ref{eq:1.1})  we consider the radiation solutions $\psi$ that is the solutions satisfying the Sommerfeld's radiation condition
\beq\label{eq:1.2}  
|x|(\frac{\partial}{\partial|x|}-i\kappa)\psi(x)\to 0\ \ \mbox{as} \ |x|\to +\infty,
\eeq
uniformly in $x/|x|$.  We assume that $\psi \in C^2$ in the closure of ${\cal U}$.

Let  
\beq\label{eq:1.3}  
 L=L_{x_0, \theta}=\{x\in\RR^3:\ \ x=x(s)=x_0+s\theta,\ \ 0<s <+\infty  \}, \ \ x_0\in\RR^3, \ \  \theta\in\mathbb S^2,
\eeq
where  $\mathbb S^2$  is the unit sphere in $\RR^3$.   

In the present work we show that, for any  radiation solution $\psi$ and any ray $L\subset ~{\cal U}$, 
the restriction of $\psi$  to $L$ is uniquely determined  by its imaginary part $\Im\psi$ on an arbitrary interval of $L$;
see Theorem 1 in Section 2.

As a corollary, we also obtain that, for any  radiation solution $\psi$ and any plane $X\subset  {\cal U}$,  the restriction of $\psi$ to $X$  is uniquely determined  by  $\Im\psi$
on an arbitrary open domain of $X$; see Theorem 2  in Section 2. 

As a further corollary,  we obtain that, for any  radiation solution $\psi$ and any plane $X\subset  {\cal U}$, the solution $\psi$ on the whole ${\cal U}$ is also uniquely determined  by  $\Im\psi$
on an arbitrary open domain of $X$; see Corollary 1  in Section 2. 

These results have holographic prototypes in the recent work \cite{N5}, where the aforementioned reconstructions of  $\psi$ are considered from the intensity $|e^{ikx}+\psi|^2$ in place of $\Im\psi$.  
Here, $e^{ikx}$ is a plane wave solution of (\ref{eq:1.1}), i.e., $k\in\RR^3,\ \ |k|=\kappa$. The results of \cite{N5} solve one of old mathematical questions of holography and 
admit straightforward applications to phaseless inverse scattering.

In the present work our studies are motivated by the Gelfand-Krein-Levitan inverse problem (from boundary values of the spectral measure in the whole space)
and by passive imaging; see, e.g.,  \cite{AHN1},  \cite{AHN2},  \cite{Ber},  \cite{BSS},  \cite{GP},  \cite{G},  \cite{MH},  \cite{WL}.

The Gelfand-Krein-Levitan problem  in question (in its fixed energy version in dimension $d=3$) consists in determining  the potential $v$ in the Schr\"odinger equation
\beq\label{eq:1.4}
-\Delta\psi(x) +v(x)\psi(x)= \kappa^2 \psi(x) +\delta(x-y),\ \ x,y\in \RR^3,\ \  \kappa>0, 
\eeq
from the imaginary part of the outgoing Green function $R_v^+(x, y, \kappa)$ for one $\kappa$ and all $x$, $y$ on some part of the boundary of a domain containing the support of $v$. 
Here,  $\delta$ is the Dirac delta function.                      
In this problem  $\Im R_v^+$ is related to the spectral measure of the Schr\"odinger operator $H=-\Delta+v$.  More precisely, $H$ 
admits the following spectral decomposition in $L^2(\RR^d)$, at least, for real-valued  compactly supported $v\in  L^{\infty}(\RR^d)$:
\beq\label{eq:1.6}
H=\int_0^{\infty}\kappa^2 d\mu_{\kappa}+\sum_{j=1}^N E_j \pi_j, \ \  d\mu_{\kappa}=\frac{2}{\pi}\Im R_v^+(\kappa)\kappa d\kappa, 
\eeq
where $d\mu_{\kappa}$ is the positive part of the spectral measure for $H$, $E_j$ are nonpositive eigenvalues for $H$ and 
$\pi_j$ are orthogonal projectors on corresponding eigenspaces, $R_v^+(\kappa)=(H-\kappa^2-i0)^{-1}$ is the limiting absorption resolvent for $H$,
whose Schwartz kernel is given by $R_v^+(x, y, \kappa)$; see, e.g.,  Lemma 14.6.1 in \cite{H}.

Note that the terminology "Gelfand-Krein-Levitan problem" is not conventional.
Our motivation for this terminology is based on the Yu.M. Berezanskii's work \cite{Ber},
which is one of the very first mathematical works on the multidimensional inverse problems for differential equations.
According to  \cite{Ber},  the problem of finding  potential $v$ in the multidimensional Schr\"odinger equation 
from boundary values of the spectral measure for some fixed boundary condition  was formulated originally 
by M.G. Krein, I.M. Gelfand and B.M. Levitan at a conference on differential equations in Moscow in 1952.
In addition,   \cite{Ber}  gives, in particular,  uniqueness theorems  on such problems including  the case of the problem of determining  $v$
from boundary values (together with some normal derivatives) of the spectral measure arising in  (\ref{eq:1.6}) for all real energies on a part of the boundary,
at least, for piecewise real-analytic $v$.

In addition,   $R_v^+(x, y, \kappa)$, for fixed $y$, can be defined as the solution $\psi$ of equation (\ref{eq:1.4}) with the radiation condition  (\ref{eq:1.2}). 
For more details, see    \cite{AHN1},  \cite{Ber}  and references therein.

Equation (\ref{eq:1.4}) at fixed $\kappa$ can be also considered as the Helmholtz equation of acoustics or electrodynamics for monochromatic waves,
where complex-valued $v(x)= v(x,\kappa)$  is related to the perturbation of the refraction index.
In particular, the aforementioned mathematical problem of recovering $v$ from boundary values of  $\Im R_v^+$ arises in the framework of
passive acoustic tomography (in ultrasonics, ocean acoustics, local helioseismology).   In these  framework
$\Im R_v^+$ is related to cross correlations of wave fields generated by random sources; see, e.g., formula (49) in \cite{G}.
For more details, see  \cite{AHN1},  \cite{AHN2},  \cite{BSS},  \cite{G},  \cite{S},  \cite{WL}   and references therein.

Note that
\beq\label{eq:1.7}
R_0^+(x, y, \kappa)=\frac{e^{i\kappa|x-y|}}{4\pi|x-y|},\ \ x,y\in \RR^3,
\eeq
where $R_0^+$ is the outgoing Green function for equation (\ref{eq:1.4}) with $v \equiv 0$.

Let 
\beq\label{eq:1.8}
R_{v,sc}^+(x, y, \kappa)=R_v^+(x, y, \kappa)-R_0^+(x, y, \kappa),\ \ x,y\in \RR^3.
\eeq

Suppose that
\beq
\label{eq:1.5}
supp\ v \subset \RR^3\setminus  \overline {{\cal U}},
\eeq
where  $\overline {{\cal U}}$ is the closure of ${\cal U}$.
Then, in view of  the definitions of $R_{v,sc}^+$,  $R_v^+$, and $R_0^+$, the function $\psi=R_{v,sc}^+(x, y, \kappa)$ is a radiation solution of equation (\ref{eq:1.1})
for each $y\in \RR^3$.

Therefore, the aforementioned results on recovering a radiation solution $\psi$ from  $\Im \psi$ give a reduction of  the Gelfand-Krein-Levitan problem
(of inverse spectral theory and passive imaging  in dimension $d=3$) to the inverse scattering problem of finding $v$ in (\ref{eq:1.4})  from boundary values of  $R_v^+$.
Note that  studies on the latter problem also go back to \cite{Ber}.

This reduction and known results imply, in particular,  that $v$ in (\ref{eq:1.4}) is uniquely determined by $\Im R_v^+(x, y, \kappa)$
for one $\kappa$ and all $x$, $y$ on an arbitrary open domain ${\cal D}$ of $X$, 
where  $supp\ v \subset  \Omega$,  $X$ is a plane in $\RR^3 \setminus  \overline { \Omega}$,  $\Omega$ is an open bounded connected domain in $\RR^3$,  $\overline { \Omega}$ is the closure of  $\Omega$;
see Theorem~3 in Section 2. 
By this result we continue studies of \cite{Ber} mentioned above and relatively recent studies of \cite{AHN1} and \cite{AHN2}.

We also consider other surfaces for measurements
instead of the planes~$X$; see Example 1 and Theorems 4 and 5 in Section 2.

The main results of this work are presented in more detail and proved in Sections 2--7.

\section{ Main results}

In this work our key result is as follows.

\begin{theorem}\label{thm:1}
Let $\psi$ be a radiation solution of equation (\ref{eq:1.1}) as in  (\ref{eq:1.2}).
Let $L$ be a ray as in (\ref{eq:1.3}) such that  $L \subset {\cal U}$, where ${\cal U}$ is the region in  (\ref{eq:1.1}). Then $\psi$ on $L$ is  uniquely determined by $\Im \psi$ on   $\Lambda$,
where   $\Lambda$  is an arbitrary non-empty open interval of~$L$.
\end{theorem}

As a corollary, we also get, in particular, the following result.

\begin{theorem}\label{thm:2}
Let $\psi$  be a radiation solution of equation (\ref{eq:1.1}) as in  (\ref{eq:1.2}).
Let $X$ be a two-dimensional plane in $\RR^3$ such that $X\subset {\cal U}$.
Then $\psi$ on $X$ is  uniquely determined by  $\Im \psi$ on  ${\cal D}$,
where  ${\cal D}$  is an arbitrary non-empty open domain of $X$.
\end{theorem}

Theorem 1 is proved in Section 4 using the Atkinson-Wilcox expansion of  \cite{A}, \cite{W} for the radiation solutions  $\psi$ of equation (\ref{eq:1.1}), and a modified version of  holographic techniques  of  \cite{N1},  \cite{N5}.
In particular, in this proof we use Proposition~1 of Section 3, which yields a two-point approximation for $\psi$
in terms of $\Im \psi$. This  two-point approximation  is also of independent interest.

Note that $\psi$ and $\Im\psi$ are real-analytic on ${\cal U}$, and, therefore, on $L$ in Theorem~1 and on $X$ in Theorem 2.
Because of this analyticity,  Theorem~1 reduces to the case when  $\Lambda=L$ and Theorem 2 reduces to the case when ${\cal D}=X$.

Using this reduction, Theorem 2 is proved as follows.
We assume that ${\cal D}=X$.  
Then to determine $\psi$ at  an arbitrary $x\in X$,  we consider a ray $L=L_{x_0, \theta}\subset X$  such that $x\in L$  and  use Theorem 1 for this $L$.

\begin{corollary}
Under the assumptions of Theorem 2,  the  imaginary part  $\Im\psi$ on  ${\cal D}$, uniquely determines $\psi$  
in the entire region ${\cal U}$.  
\end{corollary}

Corollary 1 follows from  Theorem 2, formula (\ref{eq:3.8a}) recalled in Subsection~3.3, and analyticity of  $\psi$ in ${\cal U}$.

\begin{remark}
In the one dimensional case, an analog of Theorem 1 follows from a very simple form of one dimensional  radiation solutions.
In particular, in this case an analog of the two-point approximation of Proposition 1 is exact.
The two dimensional case is considered in  \cite{NN2} using  results of  \cite{K},  \cite{NN1} in place of the three dimensional Atkinson-Wilcox expansion (\ref{eq:3.2}) used in the present work.
We expect that the case of dimension $d>3$ is similar to the three-dimensional case if $d$ is odd and is similar to the two-dimensional case if $d$ is even.
\end{remark}

Theorems 1 and 2,  and Corollary 1  have  holographic prototypes  in  \cite{N5}; see Introduction for some comments in this connection.

Using Theorem 2 and known results on direct and inverse scattering we obtain the following global uniqueness theorem for the Gelfand-Krein-Levitan inverse problem mentioned in Introduction.
\begin{theorem}\label{thm:3}
Let $v \in L^{\infty}(\RR^3)$, $supp\ v \subset  \Omega$,  and $X \subset \RR^3 \setminus  \overline { \Omega}$, where $\Omega$ is an open bounded connected domain in $\RR^3$,  $\overline { \Omega}$ is the closure of  $\Omega$,
and $X$ is a two-dimensional plane.   Let property (\ref{eq:3.11}) hold for fixed $\kappa>0$ and $R_v^+(x, y, \kappa)$ be the outgoing Green function for equation (\ref{eq:1.4}). 
Then $v$ is  uniquely determined by $\Im R_v^+(\cdot, \cdot, \kappa)$ on ${\cal D} \times {\cal D}$,
where  ${\cal D}$  is an arbitrary non-empty open domain of $X$.
\end{theorem}
In Theorem 3 we do not assume that $v$ is real-valued, 
but we assume that property (\ref{eq:3.11}) formulated in Subsection 3.4 holds.

\begin{remark}
Under the conditions that $v$ is real-valued,  $v \in L^{\infty}(\RR^d)$, $supp\ v \subseteq  \overline { \Omega}$, 
where $\Omega$ is an open bounded  domain in $\RR^d$  with $\partial\Omega \in C^{2,1}$, $\overline { \Omega}=\Omega \cup \partial\Omega$, $d \geq 2$, 
a local uniqueness for determining $v$ from
$\Im R_v^+(\cdot, \cdot, \kappa)$ on $\partial\Omega \times \partial\Omega$ for one $\kappa$
is proved in  \cite{AHN1}. The assumption that $v$ is real-valued is very essential  in this proof.
\end{remark}

\begin{remark}
In connection with applications to helioseismology, it is natural to assume that $v$ is complex-valued, $v \in L^{\infty}(\RR^3)$,  $v(x)=\tilde v(|x|)$, $\tilde v(r)=\alpha/r$, $r \geq r_0$,
for some constants $\alpha \in \RR$ and $r_0>0$.  Let $M_r=\mathbb S^2_r \times \mathbb S^2_r$ and $r_2>r_1>r_0$,
where $\mathbb S^2_r$ is defined by (\ref{eq:2.1}).
Under these assumptions,  a global uniqueness for determining $v$ from
$\Im R_v^+(\cdot, \cdot, \kappa)$ on $M_{r_1} \cup M_{r_2}$,  for fixed $\kappa$ and nonsingular pair $r_1$, $r_2$,
is proved in  \cite{AHN2}.
\end{remark}

One can see that  
Theorem 3 contains a principal progress  on the Gelfand-Krein-Levitan inverse problem in comparison with the results of  \cite{AHN1},  \cite{AHN2} mentioned in Remarks 2 and 3.
In comparison with the aforementioned result of  \cite{AHN1}, Theorem 3 is global and  does not assume that $v$ is real-valued.
In comparison with the aforementioned result  \cite{AHN2}, Theorem 3 does not assume that $v$ is spherically symmetric.

Theorem 3 is proved in Section 5.
In particular, in this proof we use Proposition 2 of Section 5,
which stays that  $v$ is  uniquely determined by $R_v^+(\cdot, \cdot, \kappa)$ on $X\times X$ for fixed  $\kappa$.
To our knowledge, Proposition 2 is slightly different from the results reported in the literature.
Therefore, just in case, for completeness of presentation this proposition is also proved in  Section 5.

\vskip 4pt
In connection with other surfaces of measurements instead of the planes $X$ our results are as follows.

The results of Theorems 1 and 2  don't hold for some other curves in place of $L$ and surfaces in place of $X$.
An example is as follows.

Let
\beq
\label{eq:2.1}
\mathbb S^2_r=\{x\in\RR^3:\ \ |x|=r\},\ \ r>0. 
\eeq
\begin{example}
Let $\psi(x)= G^+(x,\kappa)$, where $G^+$ is defined by (\ref{eq:3.8b}).  Then  $\psi$ is a non-zero radiation solution of equation   (\ref{eq:1.1})  if $\{0\} \in \RR^3\setminus  \overline {{\cal U}}$, but $\Im\psi \equiv 0$ on the  spheres $\mathbb S^2_r$ for 
$r=n\pi/ \kappa$, $n \in\NN$.
\end{example}
Nevertheless, the uniqueness results of Theorem 2, Corollary 1, and Theorem 3  remain valid for many interesting surfaces instead of the planes $X$.
Suppose that
\begin{align}
&{\cal B}\  \mbox{is\ an\ open\ bounded\ domain\ in}\ \RR^3, \label{eq:2.2}\\
&Y=\partial {\cal B}\  \mbox{is\ real-analytic\ and\ connected}. \nonumber
\end{align}

In particular, we have the following uniqueness theorems.
\begin{theorem}\label{thm:4}
Let $\psi$  be a radiation solution of equation (\ref{eq:1.1}) as in  (\ref{eq:1.2}).
Let $Y$ be a surface as in (\ref{eq:2.2}),
where $\kappa$ is not a Dirichlet eigenvalue for ${\cal B}$, and $\overline{{\cal B}}= {\cal B} \cup Y  \subset {\cal U}$.
Then $\psi$ on ${\cal U}$ is  uniquely determined by  $\Im \psi$ on  ${\cal D}$,
where  ${\cal D}$  is an arbitrary non-empty open domain of $Y$.
\end{theorem}

\begin{theorem}\label{thm:5}
Let $v \in L^{\infty}(\RR^3)$, $supp\ v \subset  \RR^3 \setminus  \overline {{\cal U}}$, where  ${\cal U}$ is as in (\ref{eq:1.1}),
and property (\ref{eq:3.11}) holds for fixed $\kappa>0$.
Let $Y$ be a surface as in (\ref{eq:2.2}),
where $\kappa$ is not a Dirichlet eigenvalue for ${\cal B}$, and $\overline{{\cal B}}= {\cal B} \cup Y  \subset {\cal U}$.
Then $v$ is  uniquely determined by $\Im R_v^+(\cdot, \cdot, \kappa)$ on ${\cal D} \times {\cal D}$,
where  ${\cal D}$  is an arbitrary non-empty open domain of $Y$.
\end{theorem}
Theorems 4 and 5 are proved in Sections 6 and 7.

Note  that the determinations in Theorems 1, 2, 3 and especially in  Theorems 4 and 5 include analytic continuations. 
Studies on more stable reconstructions appropriate for numerical  implementation  will be continued  elsewhere.
In this respect one can use, in particular, results of Subsection~3.2.

\section{Preliminaries}

\subsection{The Atkinson-Wilcox expansion }

Let
\beq
\label{eq:3.1}
B_r=\{x\in\RR^3:\ \ |x|<r\},\ \ r>0. 
\eeq
 If  $\psi$ is a radiation solution of  equation (\ref{eq:1.1}) and $\RR^3\setminus B_r  \subset {\cal U}$, then 
the following Atkinson-Wilcox expansion holds:
\begin{equation}
\label{eq:3.2}  
\psi(x)=\frac{e^{i\kappa|x|}}{|x|}\sum\limits_{j=1}^{\infty}\frac{f_j(\theta)}{|x|^{j-1}} \ \ \mbox{\rm for}\ \ x\in \RR^3\setminus B_r,\ \ \theta=\frac{x}{|x|},
\end{equation}
where the series converges absolutely and uniformly; see  \cite{A}, \cite{W}.

\subsection{A two-point approximation for $\psi$}

Let 
\begin{equation}
\label{eq:3.3}  
I(x)=|x|\Im\psi(x),\ \ x\in{\cal U},
\end{equation}
where $\psi$ is a radiation solution of equation (\ref{eq:1.1}).

We have that 
 \begin{equation}
\label{eq:3.4}  
2iI(x)=e^{i\kappa|x|}f_1(\theta)-e^{-i\kappa|x|} \overline {f_1(\theta)}+O(|x|^{-1})\ \ \mbox{\rm as} \ |x|\to +\infty,  
\end{equation}
uniformly in $\theta=x/|x|$, where $f_1$ is the leading coefficient in (\ref{eq:3.2}).

\begin{proposition}\label{prn:1}
Let $\psi$ be a radiation solution of equation (\ref{eq:1.1}). Then
 \begin{align}
&f_1(\theta)= \frac{1}{sin(\kappa\tau)}\bigl(-e^{-i\kappa|y|}I(x)+e^{-i\kappa|x|}I(y)+O(|x|^{-1})\bigr), \label{eq:3.5}\\
&x,y\in L_{x_0, \theta},\ \ x_0=0,\ \ y=x+\tau\theta, \ \ \ \  \theta\in\mathbb S^2,\ \  \tau>0, \nonumber
\end{align}
uniformly in $\theta$, where $f_1$ is the coefficient in (\ref{eq:3.4}), $I$ is defined by (\ref{eq:3.3}),  $L_{x_0, \theta}$ is the ray defined by (\ref{eq:1.3}),
and $sin(\kappa\tau)\ne 0$ for fixed $\tau$.
\end{proposition}

Formula (\ref{eq:3.5})  is a two-point approximation for $f_1$ and together with (\ref{eq:3.2}) also gives a two-point approximation for $\psi$ 
in terms of $I$. 
For phaseless inverse scattering and holography  formulas of such a type go back to \cite{N1} (see also \cite{NS2}, \cite{NSh},  \cite{N5}).

\begin{remark}
If an arbitrary function $I$ on  $L_{0, \theta}$ satisfies (\ref{eq:3.4}),
then formula (\ref{eq:3.5}) holds, for fixed $\theta\in\mathbb S^2$.
\end{remark}

We obtain (\ref{eq:3.5}) from the system of equations for $f_1$ and $\overline {f_1}$:
 \begin{align}
&e^{i\kappa|x|}f_1(\theta)-e^{-i\kappa|x|} \overline {f_1(\theta)}=2iI(x)+O(|x|^{-1}), \label{eq:3.6}\\
&e^{i\kappa|y|}f_1(\theta)-e^{-i\kappa|y|} \overline {f_1(\theta)}=2iI(y)+O(|y|^{-1}), \nonumber
\end{align}
where $x$, $y$ are as in (\ref{eq:3.5}).  In particular, we use that $|y|=|x|+\tau$ and that 
\begin{equation}
\label{eq:3.7}  
D=2i sin(\kappa\tau),
\end{equation}
where $D$ is the determinant of system (\ref{eq:3.6}).

In turn,  (\ref{eq:3.6}) follows from (\ref{eq:3.4}).

\subsection{A Green type formula}

The following formula holds:
\begin{align}
&\psi(x)=  2\int_{X}\frac{\partial G^+(x-y,\kappa)}{\partial \nu_{y}}\psi(y)dy,\ \ x\in V_X,  \label{eq:3.8a}\\
&G^+(x,\kappa)= -\frac{e^{i\kappa|x|}}{4\pi|x|},\ \ x\in \RR^3,  \label{eq:3.8b}
\end{align}
where $\psi$ is a radiation solution of equation (\ref{eq:1.1}), $X$ and $V_X$ are plane and open half-space in ${\cal U}$,
where $X$ is the boundary of $V_X$, $\nu$  is the outward normal to $X$ relative to $V_X$; see, for example, formula (5.84) in \cite{BR}.

Recall that $R_0^+(x, y, \kappa)=-G^+(x-y,\kappa)$,  where $R_0^+$ is the outgoing Green function for equation (\ref{eq:1.4}) with $v \equiv 0$.

For completeness of presentation note that formula (\ref{eq:3.8a}) follows from the formula
\begin{align}
&\psi(x)=  \int_{X}\frac{\partial G^+_X(x,y,\kappa)}{\partial \nu_{y}}\psi(y)dy,\ \ x\in V_X,  \label{eq:3.9a}\\
&G^+_X(x,y,\kappa)= G^+(x-y,\kappa)-G^+(x-y^*,\kappa),\ \ x,y\in V_X \cup X,   \nonumber
\end{align}
where $y^*$ is symmetric to $y$ with respect to $X$. 
The point is that $G^+_X$ is the  Green function for the Helmholtz  operator $\Delta+\kappa^2$ in $V_X$
with Dirichlet boundary condition on $X$ and Sommerfeld radiation condition at infinity.

\subsection{Some facts of direct scattering}

We consider equation (\ref{eq:1.4}) assuming for simplicity that
\begin{align}
&v \in L^{\infty}(\Omega),\  v \equiv 0\  \mbox{on}\  \RR^3\setminus  \overline {\Omega} \label{eq:3.9}\\
&\Omega\ \ \mbox{is\ an\ open\ bounded\ connected\ domain\ in}\ \RR^3. \nonumber
\end{align}

The outgoing Green function $R_v^+$ for equation (\ref{eq:1.4}) satisfies the integral equation
\begin{equation}
\label{eq:3.10}  
R_v^+(x, y, \kappa)=-G^+(x-y,\kappa) +\int_{\Omega}G^+(x-z,\kappa)v(z)R_v^+(z, y, \kappa)dz,
\end{equation}
where $x,y \in \RR^3$, $G^+$ is given by (\ref{eq:3.8b}).

Actually, in addition to (\ref{eq:3.9}), we assume that, for fixed  $\kappa>0$,
\begin{equation}
\label{eq:3.11}  
\mbox{equation}\ (\ref{eq:3.10})\ \mbox{is\ uniquely\ solvable\ for}\ R_v^+(\cdot , y, \kappa) \in L^2(\Omega).
\end{equation}
In particular, it is known that if $v$ satisfies (\ref{eq:3.9}) and is real-valued (or  $\Im v \leq 0$),  then (\ref{eq:3.11}) is fulfilled automatically;
see, for example, \cite{CK}.

We also consider the scattering wave functions $\psi^+$ for the homogeneous equation (\ref{eq:1.4}) (i.e.  without $\delta$):
\begin{equation}
\label{eq:3.12}  
\psi^+=\psi^+(x,\theta,\kappa)=e^{i \kappa\theta x}+\psi^+_{sc}(x,\theta,\kappa),\ \ x \in \RR^3,\ \  \theta \in \mathbb S^2,
\end{equation}
where $\psi^+_{sc}$ satisfies the radiation condition (\ref{eq:1.2}) at fixed $\theta$.

The following formulas hold:
\begin{equation}
\label{eq:3.13}  
R_v^+(x, y, \kappa)=R_v^+(y, x, \kappa),\ \ x, y\in \RR^3;
\end{equation}
\begin{equation}
\label{eq:3.14}  
R_v^+(x, y, \kappa)=-\frac{e^{i\kappa|x|}}{4\pi|x|}\psi^+(y, -\frac{x}{|x|},\kappa)+ O\bigl(\frac{1}{|x|^2}\bigr) \ \ \mbox{as} \ |x|\to +\infty
\end{equation}
at fixed $y$;
\begin{equation}
\label{eq:3.15}  
\psi^+_{sc}(x,\theta,\kappa)=\frac{e^{i\kappa|x|}}{|x|}A(\theta, \frac{x}{|x|},\kappa)+ O\bigl(\frac{1}{|x|^2}\bigr) \ \ \mbox{as} \ |x|\to +\infty
\end{equation}
at fixed $\theta$, where $A$ arising in  (\ref{eq:3.15})  is the scattering amplitude for the homogeneous equation (\ref{eq:1.4})
and  is defined on $\mathbb S^2 \times \mathbb S^2$ at fixed $\kappa$. 

In view of (\ref{eq:1.7}), (\ref{eq:1.8}), (\ref{eq:3.12})-(\ref{eq:3.14}),  we also have that
\begin{equation}
\label{eq:3.16}  
R_{v,sc}^+(x, y, \kappa)=R_{v,sc}^+(y, x, \kappa),\ \ x, y\in \RR^3;
\end{equation}
\begin{equation}
\label{eq:3.17}  
R_{v,sc}^+(x, y, \kappa)=-\frac{e^{i\kappa|x|}}{4\pi|x|}\psi_{sc}^+(y, -\frac{x}{|x|},\kappa)+ O\bigl(\frac{1}{|x|^2}\bigr) \ \ \mbox{as} \ |x|\to +\infty
\end{equation}
at fixed $y$.

In connection with aforementioned facts concerning $R_v^+$ and $\psi^+$ see, e.g., Section~1 of Chapter IV of  \cite{FM}.

\begin{remark}
It is well known that, under assumptions  (\ref{eq:3.9}),  (\ref{eq:3.11}),  the scattering amplitude $A=A(\cdot, \cdot, \kappa)$ is real analytic on $\mathbb S^2 \times \mathbb S^2$.
\end{remark}


\section{Proof of  Theorem 1}

\subsection{Case  $L\subseteq L_{0, \theta}$}
First, we give the proof for the case when $L\subseteq L_{0, \theta}$.
In this case it is essentially sufficient to prove that  $\Im\psi$ on $L$ uniquely determines  $f_j(\theta)$  in (\ref{eq:3.2}) for all $j$. 
Such a determination is presented below in this subsection.
The rest follows from the convergence of the series in (\ref{eq:3.2}) and analyticity of $\psi$ and  $\Im\psi$  on~$L$.

The determination of $f_1$ follows from (\ref{eq:3.5}).

Suppose that $f_1$,...,$f_n$ are determined, then the determination of $f_{n+1}$ is as follows.

Let
\begin{equation}
\label{eq:4.1}  
\psi_n(x)=\frac{e^{i\kappa|x|}}{|x|}\sum\limits_{j=1}^{n}\frac{f_j(\theta)}{|x|^{j-1}}, \ \ \mbox{\rm where}\ \ \theta=\frac{x}{|x|},
\end{equation}
\begin{equation}
\label{eq:4.2}  
I_n(x)=|x|\Im\psi_n,
\end{equation}
\begin{equation}
\label{eq:4.3}  
J_n(x)=|x|^n(I(x)-I_n(x)),
\end{equation}
where $x$ is as in (\ref{eq:3.2}), $I(x)$ is defined by (\ref{eq:3.3}).

We have that 
\begin{equation}
\label{eq:4.4}
2iI(x)=2iI_n(x)+  \frac{e^{i\kappa|x|}}{|x|^n}{f_{n+1}(\theta)}- \frac{e^{-i\kappa|x|}}{|x|^n}\overline{f_{n+1}(\theta)}+O(|x|^{-n-1}),
\end{equation}
\begin{equation}
\label{eq:4.5}  
2iJ_n(x)=e^{i\kappa|x|}{f_{n+1}(\theta)} - e^{-i\kappa|x|}\overline{f_{n+1}(\theta)}+O(|x|^{-1}),
\end{equation}
as $|x|\to +\infty$,  where $I$ is defined by (\ref{eq:3.3}).

Due to (\ref{eq:4.5}) and Remark 4,  we get
\begin{align}
&f_{n+1}(\theta)= \frac{1}{sin(\kappa\tau)}\bigl(-e^{-i\kappa|y|}J_n(x)+e^{-i\kappa|x|}J_n(y)+O(|x|^{-1})\bigr), \label{eq:4.6}\\
&x,y\in L_{x_0, \theta},\ \ x_0=0,\ \ y=x+\tau\theta, \ \ \ \  \theta\in\mathbb S^2,\ \  \tau>0, \nonumber
\end{align}
assuming that  $sin(\kappa\tau)\ne 0$ for fixed $\tau$
(where  the parameter $\tau$ can be always fixed in such a way for fixed $\kappa>0$).

Formulas (\ref{eq:3.3}), (\ref{eq:4.1})-(\ref{eq:4.3})  and (\ref{eq:4.6}) determine  $f_{n+1}$, give the step of induction for finding all $f_j$, and complete the proof of Theorem 1 for  the case $L\subseteq L_{0, \theta}$.

\subsection{General case}
The general case reduces to the case of Subsection 4.1 by the change of variables 
\begin{align}
&x'=x-q \label{eq:4.7}\\
&\mbox{\rm for some fixed }\ \ q\in\RR^3\ \ \mbox{\rm such that}\ \ L\subseteq L_{q, \theta}.   \nonumber
\end{align}
In the new variables $x'\in{\cal U'}={\cal U}-q$, we have that:
\begin{equation}
\label{eq:4.8}  
\psi= \frac{e^{i\kappa|x'|}}{|x'|}\sum\limits_{j=1}^{\infty}\frac{f_j'(\theta)}{|x'|^{j-1}} \ \ \mbox{\rm for}\ \ x'\in \RR^3\setminus B_{r'}, \ \  \theta=\frac{x'}{|x'|},  
\end{equation}
for some new $f_j'$,
where  $r'$ is such that $\RR^3\setminus B_{r'} \subset {\cal U'}$;
\begin{equation}
\label{eq:4.9}  
L\subseteq L_{q, \theta}=L_{0, \theta}.   
\end{equation}

In addition, the series in (\ref{eq:4.8}) converges absolutely and uniformly.

In view of (\ref{eq:4.8}), (\ref{eq:4.9}),
we complete the proof of Theorem 1 by repeating the proof of Subsection 4.1.
\begin{remark}
Our proof of Theorem 1 has a holographic prototype in \cite{N5}.
Additional formulas for finding $f_j$ from  $\Im\psi$  on $L$
can be obtained  proceeding also from  the approaches  of  \cite{N4}, \cite{NS2}, \cite{NSS}.
\end{remark}

\section{Proof of  Theorem 3}

Under our assumptions on  $v$, $\Omega$, $X$, and $\kappa$, we have, in particular, that
\begin{align}
&R_{v,sc}^+(\cdot, y, \kappa)\  \mbox{\rm is real-analytic on}\ X\ \mbox{\rm for fixed}\ y \in X, \label{eq:5.1a}  \\
&R_{v,sc}^+(x, \cdot, \kappa)\  \mbox{\rm is real-analytic on}\ X\ \mbox{\rm for fixed}\ x \in X, \label{eq:5.1b}
\end{align}
where $R_{v,sc}^+$ is defined by  (\ref{eq:1.7}),  (\ref{eq:1.8}).
Here, we use that $R_{v,sc}^+(x, y, \kappa)$ satisfies the homogeneous equation (\ref{eq:1.4}) and, therefore, is real-analytic outside of $supp\ v$, 
and that symmetry (\ref{eq:3.16})  holds.

In view of  (\ref{eq:1.7}),  (\ref{eq:1.8}),  (\ref{eq:5.1a}),  (\ref{eq:5.1b}),   Theorem 3 reduces to the case when ${\cal D}=X$.

In turn,  Theorem 3  with ${\cal D}=X$ follows from formulas (\ref{eq:1.7}),  (\ref{eq:1.8}) and from Lemma 1 and Proposition~2 given below.

\begin{lemma} Under the conditions of Theorem 3, $\Im R_{v,sc}^+(\cdot, y, \kappa)$ on $X$,  uniquely determines $R_{v,sc}^+(\cdot, y, \kappa)$ on $X$, 
where $y \in \RR^3$, and $\Im R_{v,sc}^+(\cdot, \cdot, \kappa)$ on $X\times X$  uniquely determines $R_{v,sc}^+(\cdot, \cdot, \kappa)$ on $X\times X$.
\end{lemma}
Lemma 1 follows from Theorem 2 with  $\cal U$ such that $\overline { \Omega} \subset \RR^3\setminus  \overline {{\cal U}}$, $X \subset \cal U$, and 
the property that $\psi=R_{v,sc}^+(x, y, \kappa)$ is a radiation solution of equation (\ref{eq:1.1}) for each $y\in \RR^3$.

\begin{proposition} Under the conditions of Theorem 3, $R_{v,sc}^+(\cdot, \cdot, \kappa)$ on $X\times X$ uniquely determines $v$ (for fixed  $\kappa$).
\end{proposition}

Proposition 2 is proved as follows (for example).

First,
\begin{align}
&\psi=R_{v,sc}^+(\cdot, x', \kappa)\ {\rm on}\ X\   \mbox{\rm  uniquely determines}\  \psi=R_{v,sc}^+(\cdot, x', \kappa)\  {\rm on}\ V_X \label{eq:5.2}\\
&\mbox{\rm via formula}\ (\ref{eq:3.8a}), \ \mbox{\rm for each fixed }\  x' \in X.   \nonumber
\end{align}

Second,
\begin{align}
&R_{v,sc}^+(\cdot, x', \kappa)\  {\rm on}\ V_X \cup X\ \mbox{\rm  uniquely determines}\ \psi_{sc}^+(x',\theta,\kappa)\  \mbox{\rm for}\  \theta \in \Theta_X^+ \label{eq:5.3}\\
&\mbox{\rm via formula}\ (\ref{eq:3.17}), \ \mbox{\rm for each fixed}\  x' \in X,   \nonumber
\end{align}
where
\begin{equation}
\label{eq:5.4}  
 \Theta_X^{\pm}=\{\theta \in\mathbb S^2:\  \pm\theta\nu \geq 0 \},
\end{equation}
where $\nu$  is the outward normal to $X$ relative to $V_X$ considered as the interior of $X$.

Third,
\begin{align}
&\psi=\psi_{sc}^+(\cdot, \theta,\kappa)\ {\rm on}\ X\   \mbox{\rm  uniquely determines}\  \psi=\psi_{sc}^+(\cdot, \theta,\kappa)\  {\rm on}\ V_X \label{eq:5.5}\\
&\mbox{\rm via formula}\ (\ref{eq:3.8a}).   \nonumber
\end{align}

Fourth,
\begin{align}
&\psi_{sc}^+(\cdot, \theta,\kappa)\  {\rm on}\ V_X \cup X\ \mbox{\rm  uniquely determines}\ A(\theta, \theta', \kappa)\  \mbox{\rm for}\  \theta' \in \Theta_X^- \label{eq:5.6}\\
&\mbox{\rm via formula}\  (\ref{eq:3.15}), \ \mbox{\rm for each fixed}\  \theta \in\mathbb S^2,   \nonumber
\end{align}
where $\Theta_X^-$ is defined in (\ref{eq:5.4}).

Therefore, we get that
\begin{align}
&R_{v,sc}^+(\cdot, \cdot, \kappa)\  {\rm on}\ X\times X\ \mbox{\rm  uniquely determines}\ A(\cdot, \cdot, \kappa)\  \mbox{\rm on}\  \Theta_X^+ \times \Theta_X^-\label{eq:5.7}\\
&\mbox{\rm via}\ (\ref{eq:5.2}),\  (\ref{eq:5.3}),\  (\ref{eq:5.5}),\   (\ref{eq:5.6}). \nonumber
\end{align}
In addition,
\begin{align}
&A(\cdot, \cdot, \kappa)\  \mbox{\rm on}\  \Theta_X^+ \times \Theta_X^-\  \mbox{\rm  uniquely determines}\ A(\cdot, \cdot, \kappa)\  \mbox{\rm on}\  \mathbb S^2 \times \mathbb S^2 \label{eq:5.8}\\
&\mbox{\rm by analyticity}, \nonumber
\end{align}
in view of Remark 5.

Finally, Proposition 2  follows from (\ref{eq:5.7}), (\ref{eq:5.8}) and the result of \cite{N1988} (see also \cite{N1994}, \cite{HH}, \cite{E}) that, 
under assumptions  (\ref{eq:3.9}),  (\ref{eq:3.11}), the scattering amplitude $A$ at fixed  $\kappa$ uniquely determines $v$.
In this result the assumption that $v$ is real-valued or $\Im v \leq 0$ is not essential and can be replaced by assumption~(\ref{eq:3.11}).

This completes the proofs of Proposition 2 and Theorem 3.

\section{Proof of  Theorem 4}
Under the assumptions of Theorem 4, the functions $\psi$ and $\Im \psi$ are real-analytic in ${\cal U}$, in general, and on $Y$, in particular.

Therefore, $\Im \psi$ on  ${\cal D}$  uniquely  determines $\Im \psi$ on $Y$ by analytic continuation, taking also into account that $Y$ is real-analytic and connected.

In turn, $\Im \psi$ on $Y$ uniquely determines $\Im \psi$  in ${\cal B}$ by solving the Dirichlet problem for the Helmholtz equation.

In turn, $\Im \psi$  in ${\cal B}$  uniquely determines $\Im \psi$  in ${\cal U}$ by analytic continuation.

Finally, $\Im \psi$  in ${\cal U}$ uniquely determines $\psi$  in ${\cal U}$  in view of Corollary 1 in Section 2.

This completes the proof of Theorem 4.

\section{Proof of  Theorem 5}
Recall that $\Im R_v^+(\cdot, \cdot, \kappa)$ on ${\cal D} \times {\cal D}$ uniquely  determines $\Im R_{v,sc}^+(\cdot, \cdot, \kappa)$ on ${\cal D} \times {\cal D}$,
in view of (\ref{eq:1.7}), (\ref{eq:1.8}).

Under the assumptions of Theorem 5,  we have, in particular, that
\begin{equation}
\label{eq:7.1}  
\psi=R_{v,sc}^+(\cdot, y, \kappa)\  \mbox{\rm is a radiation solution of equation (\ref{eq:1.1})},\ \  y\in \RR^3.  
\end{equation}

Therefore, due to Theorem 4, $\Im R_{v,sc}^+(\cdot, y, \kappa)$ on  ${\cal D}$  uniquely  determines $R_{v,sc}^+(\cdot, y, \kappa)$ in ${\cal U}$.

Similarly, $\Im R_{v,sc}^+(x, \cdot, \kappa)$ on  ${\cal D}$  uniquely  determines $R_{v,sc}^+(x, \cdot, \kappa)$ in ${\cal U}$,
in view of (\ref{eq:3.16}).

Therefore, $\Im R_{v,sc}^+(\cdot, \cdot, \kappa)$ on ${\cal D} \times {\cal D}$ uniquely  determines 
$R_{v,sc}^+(\cdot, \cdot, \kappa)$ on ${\cal U} \times {\cal U}$. 

\vskip 8pt
Finally, $R_{v,sc}^+(\cdot, \cdot, \kappa)$ on ${\cal U} \times {\cal U}$ uniquely determines  $v$ by different ways.

For example: $R_{v,sc}^+(\cdot, \cdot, \kappa)$ on ${\cal U} \times {\cal U}$ uniquely determines
$A(\cdot, \cdot, \kappa)$ on $\mathbb S^2 \times \mathbb S^2$, in view of (\ref{eq:3.15}), (\ref{eq:3.17});
$A$ at fixed $\kappa$ uniquely determines $v$  as recalled at the end of proof of Theorem 3.

This completes the proof of Theorem 5.
\begin{remark}
Formulas relating $R_v^+$ on $\mathbb S^2_r \times \mathbb S^2_r$  and $A$ on $\mathbb S^2 \times \mathbb S^2$  at fixed $\kappa$,
where $v(x) = 0$ for $|x|>r$, $\mathbb S^2_r$ is defined by (\ref{eq:2.1}),  were given for the first time in \cite{Ber}.
\end{remark}

\begin{remark}
The case is also of interest when the boundary $Y$ in (\ref{eq:2.2}) is not connected but consists of two disjoint connected components 
$Y_1$ and $Y_2$, where $Y_1=\partial {\cal B}_1$, $Y_2=\partial {\cal B}_2$ and  ${\cal B}_1$,  ${\cal B}_2$ are open bounded domains such that $\RR^3 \setminus {\cal U} \subset {\cal B}_1  \subset {\cal B}_2$.
In this case Theorem 4 is valid with ${\cal D}$ replaced by ${\cal D}_1 \cup {\cal D}_2$, whereas Theorems 5 is valid with ${\cal D} \times {\cal D}$ replaced by $({\cal D}_1 \cup {\cal D}_2) \times {\cal D}_1$
(as well as with ${\cal D} \times {\cal D}$ replaced by $({\cal D}_1 \cup {\cal D}_2) \times {\cal D}_2$, where ${\cal D}_1$, ${\cal D}_2$ are arbitrary non-empty open domains of $Y_1$ and $Y_2$, respectively.
It is of interest to compare the later result  with Theorem 4.3  in the recent thesis  \cite{Mu},  which gives uniqueness for monochromatic passive imaging from cross correlations on $(Y_1  \cup Y_2) \times (Y_1  \cup Y_2)$.
These comparisons may lead to further important results.
\end{remark}


\pagebreak

\noindent
Roman G. Novikov\\
CMAP, CNRS, \'Ecole polytechnique,\\
Institut Polytechnique de Paris, Palaiseau, France;\\
IEPT RAS, 117997 Moscow, Russia RAS\\
E-mail: novikov@cmap.polytechnique.fr

\end{document}